\documentclass[twocolumn]{autart}

\usepackage{graphicx}          

\usepackage{amsmath}
\usepackage{amssymb}

\def\R{\mathbb{R}}
\def\N{\mathbb{N}}

\usepackage{graphics} 
\usepackage{epsfig} 
\usepackage{epstopdf}
\usepackage{subfigure}
\usepackage{color}

\graphicspath{{./Figures/}}

\allowdisplaybreaks[4]

\newcounter{cst}

\begin{document}

\begin{frontmatter}

\title{Nonlinear Boundary Output Feedback Stabilization of Reaction-Diffusion Equations\thanksref{footnoteinfo}} 

\thanks[footnoteinfo]{Corresponding author H.~Lhachemi. The work of C. Prieur has been partially supported by MIAI@Grenoble Alpes (ANR-19-P3IA-0003)
}

\author[CS]{Hugo Lhachemi}\ead{hugo.lhachemi@centralesupelec.fr},
\author[GIPSA-lab]{Christophe Prieur}\ead{christophe.prieur@gipsa-lab.fr}, 

\address[CS]{Universit{\'e} Paris-Saclay, CNRS, CentraleSup{\'e}lec, Laboratoire des signaux et syst{\`e}mes, 91190, Gif-sur-Yvette, France}  
\address[GIPSA-lab]{Universit{\'e} Grenoble Alpes, CNRS, Grenoble-INP, GIPSA-lab, F-38000, Grenoble, France}             

\begin{keyword}                           
Reaction-diffusion PDE, Nonlinear boundary control, Nonlinear sector condition, Output feedback, Finite-dimensional controller               
\end{keyword}

\begin{abstract}                          
This paper studies the design of a finite-dimensional output feedback controller for the stabilization of a reaction-diffusion equation in the presence of a sector nonlinearity in the boundary input. Due to the input nonlinearity, classical approaches relying on the transfer of the control from the boundary into the domain with explicit occurrence of the time-derivative of the control cannot be applied. In this context, we first demonstrate using Lyapunov direct method how a finite-dimensional observer-based controller can be designed, without using the time derivative of the boundary input as an auxiliary command, in order to achieve the boundary stabilization of general 1-D reaction-diffusion equations with Robin boundary conditions and a measurement selected as a Dirichlet trace. We extend this approach to the case of a control applying at the boundary through a sector nonlinearity. We show from the derived stability conditions the existence of a size of the sector (in which the nonlinearity is confined) so that the stability of the closed-loop system is achieved when selecting the dimension of the observer to be large enough.
\end{abstract}

\end{frontmatter}

\section{Introduction}

The impact of various input nonlinearities in the control design and stability assessment of finite-dimensional systems has been intensively studied in the literature~\cite{tarbouriech2011stability,tarbouriech2011ultimate,zaccarian2011modern}.  
These include, to cite a few that are commonly encountered in practical applications~\cite{bernstein1995chronological}, saturations, deadzones, general sector nonlinearities, etc. Despite their relative simplicity, these nonlinearities, even when applied in input of a linear time invariant system, induce many challenges such as multiple equilibrium points, existence of a region of attraction, etc~\cite{campo1990robust}. The extension of these problematics to infinite dimensional systems, and particularly to the control of partial differential equations (PDEs), has attracted much attention in the recent years. Among the early contributions in this field, one can find the study of saturation mechanisms in~\cite{lasiecka2003strong,slemrod1989feedback}. More recently, Lyapunov-based stabilization of different class of PDEs, including wave and Korteweg-de Vries equations, in the presence of cone-bounded nonlinearities in the control input have been reported in~\cite{marx2017cone,marx2017global,prieur2016wave}. Model predictive control was proposed in~\cite{dubljevic2006predictive} in order to achieve the feedback stabilization of reaction-diffusion equations in the presence of constraints. Singular perturbations techniques were reported in~\cite{el2003analysis}. 

In this paper, we are concerned with the output feedback stabilization of a 1-D reaction-diffusion equation presenting a sector nonlineary in the boundary input. The developed approach relies on spectral reduction methods~\cite{russell1978controllability} that have been intensively used for the control of parabolic PDEs in a great variety of settings~\cite{coron2004global,coron2006global,lhachemi2018feedback,lhachemi2019lmi,lhachemi2019pi,lhachemi2020boundary,orlov2004robust,prieur2018feedback}.   They were in particular used in the context of the both local stabilization and estimation of region of attraction for reaction-diffusion equations in the presence of an in-domain input saturation; see~\cite{mironchenko2020local} and~\cite{lhachemi2021local} in the context of state and output feedback, respectively. It should be noted that spectral reduction methods have been very attractive for parabolic PDEs because they allow the design of finite-dimensional state-feedback, making them particularly relevant for practical applications. However, due to the distributed nature of the state, the design of an observer is generally necessary. Such observers can be designed  using, e.g., either backstepping~\cite{krstic2008boundary} of spectral methods~\cite[Sec.~4.3.2]{liu2009elementary}. Because the resulting observers mimic the dynamics of the plant, they generally take the form of a PDE, resulting in an infinite-dimensional control strategy. In order to avoid late lumping approximations required for the implementation of observers with infinite dimensional dynamics, a number of works have been devoted to the design of finite-dimensional observer-based control strategies, in particular for parabolic PDEs~\cite{balas1988finite,curtain1982finite,grune2021finite,harkort2011finite,katz2020constructive,lhachemi2020finitePI,lhachemi2020finite,sakawa1983feedback}. 

Following ideas initially introduced in~\cite{sakawa1983feedback}, the control design strategy adopted in this paper takes the form of a finite-dimensional observer coupled with a finite-dimensional state feedback. More precisely, we leverage the controller architecture reported first in~\cite{sakawa1983feedback} augmented with the LMI-based procedure suggested in~\cite{katz2020constructive}. These control design procedures were enhanced and generalized to Dirichlet and/or Neumann boundary actuation and measurement in~\cite{lhachemi2020finite}, achieving the exponential stabilization of the PDE trajectories in $H^1$ norm. In these two latter works, the boundary control input $u$ was handled using the classical procedure consisting of a change of variable that allows the transfer of the input from the boundary into the domain~\cite[Sec.~3.3]{curtain2012introduction}. However, this procedure requires the use of the time derivative $v = \dot{u}$ of the input $u$ as an auxiliary input to design the control law. When considering a nonlinearity $\varphi$ in the application of the control input $u$, such an approach fails in general. This is because the actual input applied to the system is $u_\varphi = \varphi(u)$ and its time derivative reads $v_\varphi = \dot{u}_\varphi = \varphi'(u) \dot{u}$. Therefore, since $u$ remains the actual command to be applied in input of the system, it is generally not possible to use $v_\varphi$ as an auxiliary input to perform the control design. This is because such an approach would require the integration of this latter $u$-dynamics to obtain $u$ from the knowledge of $v_\varphi$. Such a dynamics may produce trajectories that are not well defined for all $t \geq 0$ (e.g., if $\varphi'$ vanishes at certain isolated points or on a certain interval, as possible in the class of nonlinear functions $\varphi$ considered in this paper) and rise stability assessment issues for the actual input signal $u$.

Using Lyapunov direct method, the first objective of this paper is to demonstrate how a finite-dimensional observer-based controller can be designed, without using the time derivative of the boundary input as an auxiliary command, in order to achieve the output feedback boundary stabilization of a linear reaction-diffusion equation with Robin boundary conditions. More specifically, we consider the case of a control input applying through a Robin boundary condition while the measurement takes the form of  a Dirichlet trace located at the other boundary. Note that a similar setting has been achieved in~\cite{grune2021finite} by invoking small-gain arguments for PDE trajectories evaluated in $L^2$ norm (see also~\cite{katz2020delayed} for a Lyapunov-based approach but which is limited to the specific case of a bounded output operator and a Neumann control input while considering the $L^2$ norm). In contrast, the results presented in this paper rely on the use of Lyapunov functionals and provide stability estimates in both $L^2$ and $H^1$ norms. This is achieved by first designing the control strategy on the original PDE (which only involves the input $u$) while the stability analysis is performed 1) based on an homogeneous representation of the PDE that involves $v = \dot{u}$; and 2) by extending the scaling strategy employed in~\cite{lhachemi2020finite}. The second objective consists in taking advantage of the aforementioned Lyapunov functional in order to tackle the case of a boundary control input subject to a sector nonlinearity. We show the existence of a size of the sector (in which the nonlinearity is confined) so that the proposed strategies always achieve the exponential stabilization of the plant when selecting the dimension of the observer to be large enough.

The remainder of this paper is organized as follows. Notations and properties of Sturm-Liouville operators are presented in Section~\ref{sec: preliminaries}. The preliminary study consisting of the design of a finite-dimensional observer-based controller for the linear reaction-diffusion without using the time derivative of the input as an auxiliary command input is reported in Section~\ref{sec: prel study}. The extension to the case of a reaction-diffusion equation with a sector nonlinearity in the boundary input is reported in Section~\ref{sec: Design in the presence of a sector nonlinearity}. A numerical illustration is carried out in Section~\ref{sec: numerical illustration}. Finally, concluding remarks are formulated in Section~\ref{sec: conclusion}.

\section{Notation and properties}\label{sec: preliminaries}

Spaces $\R^n$ are equipped with the usual Euclidean norm denoted by $\Vert\cdot\Vert$. The associated induced norms of matrices are also denoted by $\Vert\cdot\Vert$. For any two vectors $X$ and $Y$ of arbitrary dimensions, $ \mathrm{col} (X,Y)$ denotes the vector $[X^\top,Y^\top]^\top$. The space $L^2(0,1)$ stands for the space of square integrable functions on $(0,1)$ and is endowed with the usual inner product $\langle f , g \rangle = \int_0^1 f(x) g(x) \,\mathrm{d}x$. The associated norm is denoted by $\Vert \cdot \Vert_{L^2}$. For an integer $m \geq 1$, the $m$-order Sobolev space is denoted by $H^m(0,1)$ and is endowed with its usual norm $\Vert \cdot \Vert_{H^m}$. For a symmetric matrix $P \in\R^{n \times n}$, $P \succeq 0$ (resp. $P \succ 0$) means that $P$ is positive semi-definite (resp. positive definite).

Let $\theta_1\in(0,\pi/2]$, $\theta_2\in[0,\pi/2]$, $p \in \mathcal{C}^1([0,1])$ and $q \in \mathcal{C}^0([0,1])$ with $p > 0$ and $q \geq 0$. Let the Sturm-Liouville operator $\mathcal{A} : D(\mathcal{A}) \subset L^2(0,1) \rightarrow L^2(0,1)$ be defined by $\mathcal{A}f = - (pf')' + q f$ on the domain $D(\mathcal{A}) = \{ f \in H^2(0,1) \,:\, \cos(\theta_1) f(0) - \sin(\theta_1) f'(0) = \cos(\theta_2) f(1) + \sin(\theta_2) f'(1) = 0 \}$. The eigenvalues $\lambda_n$, $n \geq 1$, of $\mathcal{A}$ are simple, non negative, and form an increasing sequence with $\lambda_n \rightarrow + \infty$ as $n \rightarrow + \infty$. Moreover, the associated unit eigenvectors $\phi_n \in L^2(0,1)$ form a Hilbert basis. The domain of the operator $\mathcal{A}$ is equivalently characterized by $D(\mathcal{A}) = \{ f \in L^2(0,1) \,:\, \sum_{n\geq 1} \vert \lambda_n \vert ^2 \vert \langle f , \phi_n \rangle \vert^2 < +\infty \}$. Introducing $p_*,p^*,q^* \in \R$ such that $0 < p_* \leq p(x) \leq p^*$ and $0 \leq q(x) \leq q^*$ for all $x \in [0,1]$, then it holds 
$
0 \leq \pi^2 (n-1)^2 p_* \leq \lambda_n \leq \pi^2 n^2 p^* + q^*
$
for all $n \geq 1$~\cite{orlov2017general}. Moreover if $p \in \mathcal{C}^2([0,1])$, we have (see, e.g., \cite{orlov2017general,tucsnak2009observation}) that $\phi_n (\xi) = O(1)$ and $\phi_n' (\xi) = O(\sqrt{\lambda_n})$ as $n \rightarrow + \infty$ for any given $\xi \in [0,1]$. Assuming further that $q > 0$, an integration by parts and the continuous embedding $H^1(0,1) \subset L^\infty(0,1)$ show the existence of constants $C_1,C_2 > 0$ such that
\begin{align}
C_1 \Vert f \Vert_{H^1}^2 \leq 
\sum_{n \geq 1} \lambda_n \langle f , \phi_n \rangle^2
= \langle \mathcal{A}f , f \rangle
\leq C_2 \Vert f \Vert_{H^1}^2 \label{eq: inner product Af and f}
\end{align}
for any $f \in D(\mathcal{A})$. The latter inequalities and the Riesz-spectral nature of $\mathcal{A}$ imply that the series expansion $f = \sum_{n \geq 1} \langle f , \phi_n \rangle \phi_n$ holds in $H^2(0,1)$ norm for any $f \in D(\mathcal{A})$. Due to the continuous embedding $H^1(0,1) \subset L^{\infty}(0,1)$, we obtain that $f(0) = \sum_{n \geq 1} \langle f , \phi_n \rangle \phi_n(0)$. We finally define, for any integer $N \geq 1$, $\mathcal{R}_N f = \sum_{n \geq N+1} \langle f , \phi_n \rangle \phi_n$.

\section{Design for linear reaction-diffusion equation}\label{sec: prel study}

Consider the reaction-diffusion system described by
\begin{subequations}\label{eq: PDE linear}
\begin{align}
& z_t(t,x) = (p(x) z_x(t,x))_x - \tilde{q}(x) z(t,x) \\
& \cos(\theta_1) z(t,0) - \sin(\theta_1) z_x(t,0) = 0 \\
& \cos(\theta_2) z(t,1) + \sin(\theta_2) z_x(t,1) = u(t) \\
& y (t) = z(t,0) \label{eq: PDE linear - Dirichlet measurement} \\
& z(0,x) = z_0(x)
\end{align}
\end{subequations}
for $t > 0$ and $x \in (0,1)$ where $\theta_1 \in (0,\pi/2]$, $\theta_2 \in [0,\pi/2]$, $p \in\mathcal{C}^2([0,1])$ with $p > 0$, and $\tilde{q} \in\mathcal{C}^0([0,1])$. Here $z(t,\cdot)$ represents the state at time $t$, $u(t)$ is the command, $y(t)$ is the measurement, and $z_0$ is the initial condition. Without loss of generality, let $q \in\mathcal{C}^0([0,1])$ and $q_c \in\R$ be such that
\begin{equation}\label{eq: writting of tilde_q}
\tilde{q}(x) = q(x) - q_c , \quad q(x) > 0  .
\end{equation}

\subsection{Spectral reduction}
Introducing the change of variable
\begin{equation}\label{eq: change of variable}
w(t,x) = z(t,x) - \frac{x^2}{\cos(\theta_2) + 2 \sin(\theta_2)} u(t) 
\end{equation}
and defining $v(t) = \dot{u}(t)$ we infer that
\begin{subequations}\label{eq: PDE linear - homogeneous}
\begin{align}
& \dot{u}(t) = v(t) \\
& w_t(t,x) = (p(x) w_x(t,x))_x - \tilde{q}(x) w(t,x) \\
& \phantom{w_t(t,x) =} \;  + a(x) u(t) + b(x) v(t) \\ 
& \cos(\theta_1) w(t,0) - \sin(\theta_1) w_x(t,0) = 0 \\
& \cos(\theta_2) w(t,1) + \sin(\theta_2) w_x(t,1) = 0 \\
& y(t) = w(t,0) \\
& w(0,x) = w_0(x)
\end{align}
\end{subequations}
where $a(x) = \frac{1}{\cos(\theta_2) + 2 \sin(\theta_2)} \{ 2p(x) + 2xp'(x) - x^2 \tilde{q}(x) \}$, $b(x) = -\frac{x^2}{\cos(\theta_2) + 2 \sin(\theta_2)}$, and $w_0(x) = z_0(x) - \frac{x^2}{\cos(\theta_2) + 2 \sin(\theta_2)} u(0)$. We define $z_n(t) = \langle z(t,\cdot) , \phi_n \rangle$, $w_n(t) = \langle w(t,\cdot) , \phi_n \rangle$, $a_n = \langle a , \phi_n \rangle$, and $b_n = \langle b , \phi_n \rangle$. In particular, one has
\begin{equation}\label{eq: link z_n and w_n}
w_n(t) = z_n(t) + b_n u(t), \quad n \geq 1 .
\end{equation}
The projection of (\ref{eq: PDE linear}) into the Hilbert basis $(\phi_n)_{n \geq 1}$ gives
\begin{equation}\label{eq: dynamics z_n}
\dot{z}_n(t) = (-\lambda_n + q_c) z_n(t) + \beta_n u(t)
\end{equation}
with 
\begin{align}
\beta_n & = a_n + (-\lambda_n+q_c)b_n \nonumber \\
& = p(1) \{ - \cos(\theta_2) \phi_n'(1) + \sin(\theta_2) \phi_n(1) \} = O(\sqrt{\lambda_n}) \label{eq: def beta_n}
\end{align}
while the projection of (\ref{eq: PDE linear - homogeneous}) reads
\begin{subequations}\label{eq: dynamics w_n}
\begin{align}
\dot{u}(t) & = v(t) \\
\dot{w}_n(t) & = (-\lambda_n + q_c) w_n(t) + a_n u(t) + b_n v(t) ,  \\
y(t) & = \sum_{n \geq 1} w_n(t) \phi_n(0) \label{eq: dynamics w_n - y}
\end{align}
\end{subequations}

\begin{rem}\label{rem: technical challenges}
Representation (\ref{eq: dynamics z_n}) is more convenient for control design since only the input $u$ appears in the dynamics. However, Lyapunov stability analysis based on this representation is difficult because $\beta_n = O(\sqrt{\lambda_n})$. Conversely, representation (\ref{eq: dynamics w_n}) is less natural for control design since both input $u$ and its time derivative $v = \dot{u}$ appear in the dynamics. Nevertheless, this representation is easier to handle in the context of a Lyapunov stability analysis because $a_n,b_n\in\ell^2(\N)$; see~\cite{katz2020constructive,lhachemi2020finite} where $v = \dot{u}$ was used as the input for control design. In this section, we demonstrate for the general setting of (\ref{eq: PDE linear}) how to perform the control design directly with $u$ as the input, based on representation (\ref{eq: dynamics z_n}), while carrying out the Lyapunov stability analysis using representation (\ref{eq: dynamics w_n}) in order to obtain stability estimates in both $L^2$ and $H^1$ norms. 
\end{rem}

\subsection{Control strategy}

Let $\delta > 0$ and $N_0 \geq 1$ be such that $-\lambda_n + q_c < - \delta < 0$ for all $n \geq N_0 + 1$. For an arbitrarily given $N \geq N_0 + 1$, we design an observer to estimate the $N$ first modes of the plant in $z$-coordinates while the state-feedback is computed based on the estimation of the $N_0$ first modes of the plant. More precisely and inspired by the early work~\cite{sakawa1983feedback}, the control strategy investigated in this section takes the form: 
\begin{subequations}\label{eq: controller}
\begin{align}
\hat{w}_n(t) & = \hat{z}_n(t) + b_n u(t) \label{eq: controller 1} \\
\dot{\hat{z}}_n(t) & = (-\lambda_n+q_c) \hat{z}_n(t) + \beta_n u(t) \label{eq: controller 2} \\
& \phantom{=}\; - l_n \left\{ \sum_{k = 1}^N \hat{w}_k(t) \phi_k(0) - y(t) \right\}  ,\; 1 \leq n \leq N_0 \nonumber \\
\dot{\hat{z}}_n(t) & = (-\lambda_n+q_c) \hat{z}_n(t) + \beta_n u(t) ,\; N_0+1 \leq n \leq N \label{eq: controller 3} \\
u(t) & = \sum_{k = 1}^{N_0} k_k \hat{z}_k(t) \label{eq: controller 4}
\end{align}
\end{subequations}
Here $l_n \in\R$ and $k_k \in\R$ are the observer and feedback gains, respectively. Compared to~\cite{lhachemi2020finite} where the observer was designed to estimate the modes $w_n$ in homogeneous coordinates (\ref{eq: PDE linear - homogeneous}), the controller architecture (\ref{eq: controller}) proposed in this paper differs as it manages to directly estimate the modes the $z_n$ of the plant in original non-homogeneous coordinates (\ref{eq: PDE linear}). This change of the structure of the observer is key to perform the control design directly with $u$ instead of $v = \dot{u}$. Note that the series expansion (\ref{eq: dynamics w_n - y}) of the system output $y(t)$ holds in $w$ coordinates but not in original $z$ coordinates. This is why the estimated system output used to compensate the error of observation in (\ref{eq: controller 2}) takes the form $\sum_{k = 1}^N \hat{w}_k(t) \phi_k(0)$. Here the terms $\hat{w}_n$ stand for the estimates of the modes in $w$ coordinates. They are obtained from the estimates $\hat{z}_n$ of the modes in $z$ coordinates through (\ref{eq: controller 1}) which mimics (\ref{eq: link z_n and w_n}).

\begin{rem}
We denote by $\hat{z}(t) \in\R^N$ the state of the observer. The well-posedness of the closed-loop system composed of (\ref{eq: PDE linear - homogeneous}) and (\ref{eq: controller}) in terms of classical solutions for initial conditions $w_0 \in D(\mathcal{A})$ and $\hat{z}(0) \in\R^N$ and defined for all $t \geq 0$, namely $(w,\hat{z}) \in \mathcal{C}^0(\R_{\geq 0};L^2(0,1) \times \R^N) \cap \mathcal{C}^1(\R_{>0};L^2(0,1) \times \R^N)$ with $w(t,\cdot) \in D(\mathcal{A})$ for all $t > 0$, is a direct consequence of~\cite[Thm.~6.3.1 and~6.3.3]{pazy2012semigroups}. Moreover, from the proof of ~\cite[Thm.~6.3.1]{pazy2012semigroups}, we have $\mathcal{A} w \in \mathcal{C}^0(\R_{>0};L^2(0,1))$ and $\mathcal{A}^{1/2}w \in \mathcal{C}^0(\R_{\geq 0};L^2(0,1))$.
\end{rem}

\subsection{Model for stability analysis}\label{subsec: Model for stability analysis}

We define $e_n = z_n - \hat{z}_n$ for all $1 \leq n \leq N$. From (\ref{eq: controller 1}-\ref{eq: controller 2}) and using (\ref{eq: link z_n and w_n}) and (\ref{eq: dynamics w_n - y}), we infer that
\begin{equation}\label{eq: controller 2 bis}
\dot{\hat{z}}_n 
= (-\lambda_n + q_c) \hat{z}_n + \beta_n u + l_n \sum_{k=1}^{N} \phi_k(0) e_k + l_n \zeta
\end{equation}
for $1 \leq n \leq N_0$ where $\zeta(t) = \sum_{n \geq N+1} w_n(t) \phi_n(0)$. We define first the scaled quantities $\tilde{z}_n = \hat{z}_n / \lambda_n$ and, as in~\cite{lhachemi2020finite}, $\tilde{e}_n = \sqrt{\lambda_n} e_n$. We then introduce $\hat{Z}^{N_0} = \begin{bmatrix} \hat{z}_1 & \ldots & \hat{z}_{N_0} \end{bmatrix}^\top$, $E^{N_0} = \begin{bmatrix} e_1 & \ldots & e_{N_0} \end{bmatrix}^\top$, $\tilde{Z}^{N-N_0} = \begin{bmatrix} \tilde{z}_{N_0 + 1} & \ldots & \tilde{z}_{N} \end{bmatrix}^\top$, and $\tilde{E}^{N_0} = \begin{bmatrix} \tilde{e}_{N_0 +1} & \ldots & \tilde{e}_{N} \end{bmatrix}^\top$. We obtain from (\ref{eq: controller 4}) that
\begin{equation}\label{eq: command input - matrix form}
u = K \hat{Z}^{N_0}
\end{equation}
where $K = \begin{bmatrix} k_1 & \ldots & k_{N_0} \end{bmatrix}$. Next, we infer from (\ref{eq: controller}) and (\ref{eq: controller 2 bis}) that 
\begin{subequations}\label{eq: truncated model - 4 ODEs}
\begin{align}
\dot{\hat{Z}}^{N_0} & = (A_0 + \mathfrak{B}_0 K) \hat{Z}^{N_0} + LC_0 E^{N_0} \label{eq: truncated model - 4 ODEs - 1} \\
& \phantom{=}\; + L\tilde{C}_1 \tilde{E}^{N-N_0} + L \zeta \nonumber \\
\dot{E}^{N_0} & = ( A_0 - L C_0 ) E^{N_0} - L \tilde{C}_1 \tilde{E}^{N-N_0} - L \zeta \\
\dot{\tilde{Z}}^{N-N_0} & = A_1 \tilde{Z}^{N-N_0} + \tilde{\mathfrak{B}}_1 K \hat{Z}^{N_0} \\
\dot{\tilde{E}}^{N-N_0} & = A_1 \tilde{E}^{N-N_0}
\end{align}
\end{subequations}
where the different matrices are defined by $A_0 = \mathrm{diag}(-\lambda_1 + q_c , \ldots , -\lambda_{N_0} + q_c)$, $A_1 = \mathrm{diag}(-\lambda_{N_0+1} + q_c , \ldots , -\lambda_{N} + q_c)$, $\mathfrak{B}_0 = \begin{bmatrix} \beta_1 & \ldots & \beta_{N_0} \end{bmatrix}^\top$, $\tilde{\mathfrak{B}}_1 = \begin{bmatrix} \frac{\beta_{N_0 +1}}{\lambda_{N_0 +1}} & \ldots & \frac{\beta_N}{\lambda_N} \end{bmatrix}^\top$, $C_0 = \begin{bmatrix} \phi_1(0) & \ldots & \phi_{N_0}(0) \end{bmatrix}$, $\tilde{C}_1 = \begin{bmatrix} \frac{\phi_{N_0 +1}(0)}{\sqrt{\lambda_{N_0 +1}}} & \ldots & \frac{\phi_{N}(0)}{\sqrt{\lambda_{N}}} \end{bmatrix}$, and $L = \begin{bmatrix} l_1 & \ldots & l_{N_0} \end{bmatrix}^\top$. Therefore, defining the vector
\begin{equation}\label{eq: truncated model - def X}
X = \mathrm{col}\left( \hat{Z}^{N_0} , E^{N_0} , \tilde{Z}^{N-N_0} , \tilde{E}^{N-N_0} \right) ,
\end{equation}
we infer that
\begin{equation}\label{eq: truncated model}
\dot{X} = F X + \mathcal{L} \zeta
\end{equation}
where
\begin{equation*}
F =
\begin{bmatrix}
A_0 + \mathfrak{B}_0 K & LC_0 & 0 & L\tilde{C}_1 \\
0 & A_0 - L C_0 & 0 & - L\tilde{C}_1 \\
\tilde{\mathfrak{B}}_1 K & 0 & A_1 & 0 \\
0 & 0 & 0 & A_1
\end{bmatrix} ,\;
\mathcal{L} =
\begin{bmatrix}
L \\ -L \\ 0 \\ 0
\end{bmatrix} .
\end{equation*}

\begin{rem}
The pairs $(A_0,\mathfrak{B}_0)$ and $(A_0,C_0)$ satisfy the Kalman condition. Indeed, since $A_0$ is diagonal with simple eigenvalues, the Kalman conditions hold if and only if, from the definition of the matrices $\mathfrak{B}_0$ and $C_0$, $\beta_n = p(1) \{ - \cos(\theta_2) \phi_n'(1) + \sin(\theta_2) \phi_n(1) \} \neq 0$ and $\phi_n(0) \neq 0$ for all $1 \leq n \leq N_0$. From the definition of the eigenvectors $\phi_n$ and by Cauchy uniqueness, this is indeed case. Hence, we can always compute a feedback gain $K\in\R^{1 \times N_0}$ and an observer gain $L\in\R^{N_0}$ such that $A_0 + \mathfrak{B}_0 K$ and $A_0 - L C_0$ are Hurwitz with eigenvalues that have a real part strictly less than $-\delta<0$. In that case, from the above definition of the matrix $F$, one can observe that matrix $F$ is Hurwitz with eigenvalues that have a real part strictly less than $-\delta<0$.
\end{rem}

Finally, defining $\tilde{X} = \mathrm{col}\left( X , \zeta \right)$ and based on (\ref{eq: command input - matrix form}) and (\ref{eq: truncated model - 4 ODEs - 1}), we also have 
\begin{equation}\label{eq: derivative v of command input u}
u = \tilde{K} X  , \quad  \quad v = \dot{u} = K \dot{\hat{Z}}^{N_0} = E \tilde{X}
\end{equation}
with $E = K \begin{bmatrix} A_0 + \mathfrak{B}_0 K & LC_0 & 0 & L\tilde{C}_1 & L \end{bmatrix}$ and $\tilde{K} = \begin{bmatrix} K & 0 & 0 & 0 \end{bmatrix}$.

\subsection{Main stability results}

\begin{thm}\label{thm1}
Let $\theta_1 \in (0,\pi/2]$, $\theta_2 \in [0,\pi/2]$, $p \in\mathcal{C}^2([0,1])$ with $p > 0$, and $\tilde{q} \in\mathcal{C}^0([0,1])$. Let $q \in\mathcal{C}^0([0,1])$ and $q_c \in\R$ be such that (\ref{eq: writting of tilde_q}) holds. Let $\delta > 0$ and $N_0 \geq 1$ be such that $-\lambda_n + q_c < - \delta$ for all $n \geq N_0 + 1$. Let $K\in\R^{1 \times N_0}$ and $L\in\R^{N_0}$ be such that $A_0 + \mathfrak{B}_0 K$ and $A_0 - L C_0$ are Hurwitz with eigenvalues that have a real part strictly less than $-\delta<0$. For a given $N \geq N_0 +1$, assume that there exist $P \succ 0$, $\alpha>1$, and $\beta,\gamma > 0$ such that 
\begin{equation}\label{eq: thm1 - constraints}
\Theta_1 \preceq 0 ,\quad \Theta_2 \leq 0
\end{equation}
where
\begin{subequations}\label{eq: def theta linear}
\begin{align}
\Theta_1 & = \begin{bmatrix} F^\top P + P F + 2 \delta P + \alpha\gamma \Vert \mathcal{R}_N a \Vert_{L^2}^2 \tilde{K}^\top \tilde{K} & P\mathcal{L} \\ \mathcal{L}^\top P & - \beta \end{bmatrix} \nonumber \\
& \phantom{=}\; + \alpha\gamma \Vert \mathcal{R}_N b \Vert_{L^2}^2 E^\top E \label{eq: def theta1 linear} \\
\Theta_2 & = 2\gamma\left\{ - \left( 1 - \frac{1}{\alpha} \right) \lambda_{N+1}+q_c + \delta \right\} + \beta M_\phi
\end{align}
\end{subequations}
and with $M_\phi = \sum_{n \geq N+1} \frac{\vert \phi_n(0) \vert^2}{\lambda_n} < +\infty$. Then there exists a constant $M > 0$ such that for any initial conditions $z_0 \in H^2(0,1)$ and $\hat{z}_n(0)\in\R$ such that $\cos(\theta_1) z_0(0) - \sin(\theta_1) z_0'(0) = 0$ and $\cos(\theta_2) z_0(1) + \sin(\theta_2) z_0'(1) = K \hat{Z}^{N_0}(0)$, the trajectories of the closed-loop system composed of the plant (\ref{eq: PDE linear}) and the controller (\ref{eq: controller}) satisfy 
$$\Vert z(t,\cdot) \Vert_{H^1}^2 + \sum_{n = 1}^{N} \hat{z}_n(t)^2 \leq M e^{-2\delta t} \left( \Vert z_0 \Vert_{H^1}^2 + \sum_{n = 1}^{N} \hat{z}_n(0)^2 \right)$$ 
for all $t \geq 0$. Moreover, the constraints (\ref{eq: thm1 - constraints}) are always feasible for $N$ selected large enough.
\end{thm}

\begin{pf} 
Let the Lyapunov function candidate
\begin{equation}\label{eq: Lyapunov functional - H1 norm}
V(X,w) = X^\top P X + \gamma \sum_{n \geq N+1} \lambda_n \langle w , \phi_n \rangle^2
\end{equation}
for $X \in\R^{2N}$ and $w \in D(\mathcal{A})$. Note that the first term of $V$ accounts for the $N$ first modes of the PDE expressed in $z$-coordinates (\ref{eq: PDE linear}) while the second term accounts for the modes labeled $n \geq N+1$ of the PDE expressed in $w$-coordinates (\ref{eq: PDE linear - homogeneous}). Compared to \cite{lhachemi2020finite} where the considered Lyapunov functional only captures the dynamics of the modes in homogeneous coordinates, the rationale for considering both non-homogeneous and homogeneous representations of the plant for control design and stability analysis has been described in Remark~\ref{rem: technical challenges}. The computation of the time derivative of $V$ along classical solutions of the system composed of (\ref{eq: dynamics w_n}) and (\ref{eq: truncated model}), whose existence is provided by~\cite[Thm.~6.3.1]{pazy2012semigroups}, gives
\begin{align*}
\dot{V} & =
\tilde{X}^\top
\begin{bmatrix} F^\top P + P F & P \mathcal{L} \\ \mathcal{L}^\top P & 0 \end{bmatrix}
\tilde{X} \\
& \phantom{=}\; + 2 \gamma \sum_{n \geq N+1} \lambda_n \left\{ (-\lambda_n + q_c) w_n + a_n u + b_n v \right\} w_n .
\end{align*}
where $\tilde{X} = \mathrm{col}\left( X , \zeta \right)$. Using Young inequality, we infer for any $\alpha > 0$ that 
\begin{align*}
2\sum_{n \geq N+1} \lambda_n a_n u w_n & \leq \frac{1}{\alpha} \sum_{n \geq N+1} \lambda_n^2 w_n^2 + \alpha \Vert \mathcal{R}_N a \Vert_{L^2}^2 u^2 , \\
2\sum_{n \geq N+1} \lambda_n b_n v w_n & \leq \frac{1}{\alpha} \sum_{n \geq N+1} \lambda_n^2 w_n^2 + \alpha \Vert \mathcal{R}_N b\Vert_{L^2}^2 v^2 
\end{align*}
where, using (\ref{eq: derivative v of command input u}), $u^2 = X^\top \tilde{K}^\top \tilde{K} X$ and $v^2 = \tilde{X}^\top E^\top E \tilde{X}$. Moreover, since $\zeta = \sum_{n \geq N+1} w_n \phi_n(0)$, Cauchy-Schwarz inequality gives $\zeta^2 \leq M_\phi \sum_{n \geq N+1} \lambda_n w_n^2$. Combining the above estimates, we infer that 
$$\dot{V} + 2 \delta V \leq \tilde{X}^\top \Theta_1 \tilde{X} + \sum_{n\geq N+1}\lambda_n \Gamma_n w_n^2$$ 
where $\Gamma_n = 2\gamma\left\{ - \left( 1 - \frac{1}{\alpha} \right) \lambda_{n}+q_c + \delta \right\} + \beta M_\phi$. Recalling that $\alpha > 1$, we have $\Gamma_n \leq \Theta_2 \leq 0$ for all $n \geq N+1$. Since $\Theta_1 \preceq 0$, we have that $\dot{V} + 2 \delta V \leq 0$. Using (\ref{eq: inner product Af and f}), (\ref{eq: change of variable}), (\ref{eq: link z_n and w_n}), and (\ref{eq: command input - matrix form}), we obtain the claimed stability estimate.

It remains to show that the constraints $\Theta_1 \preceq 0$ and $\Theta_2 \leq 0$ are always feasible for $N \geq N_0 + 1$ selected large enough. To do so, we apply Lemma~\ref{lem: useful lemma} reported in appendix to the matrix $F + \delta I$. This is possible because (i) $A_0 + \mathfrak{B}_0 K + \delta I$ and $A_0 - L C_0 + \delta I$ are Hurwitz; (ii) $\Vert e^{(A_1 + \delta I)t} \Vert \leq e^{- \kappa_0 t}$ for all $t \ge 0$ with $\kappa_0 = \lambda_{N_0+1} - q_c - \delta > 0$ defined independently of $N$; and (iii) $\Vert L\tilde{C}_1 \Vert \leq \Vert L \Vert \Vert \tilde{C}_1 \Vert$ and $\Vert \tilde{\mathfrak{B}}_1 K \Vert \leq \Vert \tilde{\mathfrak{B}}_1 \Vert \Vert K \Vert$ where $K$ and $L$ are independent of the number of observed modes $N$ while $\Vert \tilde{C}_1 \Vert = O(1)$ and $\Vert \tilde{\mathfrak{B}}_1 \Vert = O(1)$ when $N \rightarrow +\infty$. Hence the solution $P \succ 0$ to the Lyapunov equation $F^\top P + P F + 2 \delta P = -I$ is such that $\Vert P \Vert = O(1)$ as $N \rightarrow + \infty$. With this choice of matrix $P$, the constraint $\Theta_1 \preceq 0$ becomes equivalent to $\Theta_{1p} + \alpha\gamma \Vert \mathcal{R}_N b \Vert_{L^2}^2 E^\top E \preceq 0$ where
\begin{align*}
\Theta_{1p} = \begin{bmatrix} - I + \alpha\gamma \Vert \mathcal{R}_N a \Vert_{L^2}^2 \tilde{K}^\top \tilde{K} & P\mathcal{L} \\ \mathcal{L}^\top P & - \beta \end{bmatrix} .
\end{align*}
We note that $\Vert \tilde{K} \Vert = \Vert K \Vert$ and $\Vert \mathcal{L} \Vert = \sqrt{2} \Vert L \Vert$ are independent of $N$. Hence, fixing the value of $\alpha > 1$ and selecting $\beta = \sqrt{N}$ and $\gamma = 1/N$ with $N \geq N_0 +1$ large enough, we infer that (i) $\Theta_2\leq 0$  and (ii) by invoking Schur complement, $\Theta_{1p} \preceq -\frac{1}{2} I$. Noting from (\ref{eq: derivative v of command input u}) that $\Vert E \Vert = O(1)$ as $N \rightarrow + \infty$, this implies that $\Theta_1 \preceq 0$ for $N \geq N_0 +1$ large enough. We have shown that the constraints (\ref{eq: thm1 - constraints}) are feasible when selecting $N \geq N_0 + 1$ to be large enough. This completes the proof.
\end{pf}

\begin{rem}\label{rem: obtention of LMIs}
For a given number of modes of the observer $N \geq N_0 +1$, the constraints (\ref{eq: thm1 - constraints}) are nonlinear w.r.t. the decision variables $P \succ 0$, $\alpha>1$, and $\beta,\gamma > 0$. However, fixing the value of $\alpha > 1$, the constraints (\ref{eq: thm1 - constraints}) now take the form of LMIs of the decision variables $P \succ 0$ and $\beta,\gamma > 0$. This latter LMI formulation of the constraints remains feasible for $N \geq N_0 +1$ selected large enough as shown in the proof of Theorem~\ref{thm1}. Note however that this approach requires to fix \emph{a priori} the value of $\alpha > 1$. Hence, the obtained value $N$ for the dimension of the observer depends in general of the fixed value of $\alpha$, which may introduce some conservatism. A second approach, which has the benefit of not introducing any conservatism compared to the original constraints (\ref{eq: thm1 - constraints}), goes as follows. Proceeding with the substitutions $P \leftarrow \gamma P$ and $\beta \leftarrow \gamma \beta$ into (\ref{eq: def theta linear}) and factoring out $\gamma > 0$, we observe that (\ref{eq: thm1 - constraints}) holds if and only if $\Theta_{1,\gamma=1} \preceq 0$ and $\Theta_{2,\gamma=1} \leq 0$ where $\Theta_{1,\gamma=1}$ and $\Theta_{2,\gamma=1}$ are obtained by setting $\gamma = 1$ into (\ref{eq: def theta linear}). In this case, $\Theta_{1,\gamma=1}$ takes the form of a LMI of the decision variables $P \succ 0$, $\alpha>1$, and $\beta > 0$ while the use of the Schur complement shows that $\Theta_{2,\gamma=1}$ can be equivalently recast into the LMI formulation:
\begin{equation*}
\begin{bmatrix}
2 \{ -\lambda_{N+1} + q_c + \delta \} + \beta M_\phi & \sqrt{2\lambda_{N+1}} \\
\sqrt{2\lambda_{N+1}} & -\alpha
\end{bmatrix}
\preceq 0 .
\end{equation*}
This LMI formulation of the problem is equivalent to (\ref{eq: thm1 - constraints}), hence is always feasible for $N \geq N_0+1$ selected to be large enough. A similar remark applies to the next theorems.
\end{rem}

\begin{rem}
The result of Theorem~\ref{thm1} in the case of the Dirichlet measurement (\ref{eq: PDE linear - Dirichlet measurement}) can easily be adapted to the case of the Neumann measurement $y(t) = z_x(t,0)$ with $\theta_1 \in [0,\pi/2)$ and $\theta_2 \in [0,\pi/2]$ by combining the approach developed in this paper along with the rescaling procedure reported in~\cite{lhachemi2020finite}.
\end{rem}

We also state below a $L^2$ version of the stability result presented in Theorem~\ref{thm1}.

\begin{thm}\label{thm2}
Let $\theta_1 \in (0,\pi/2]$, $\theta_2 \in [0,\pi/2]$, $p \in\mathcal{C}^2([0,1])$ with $p > 0$, and $\tilde{q} \in\mathcal{C}^0([0,1])$. Let $q \in\mathcal{C}^0([0,1])$ and $q_c \in\R$ be such that (\ref{eq: writting of tilde_q}) holds. Let $\delta > 0$ and $N_0 \geq 1$ be such that $-\lambda_n + q_c < - \delta$ for all $n \geq N_0 + 1$. Let $K\in\R^{1 \times N_0}$ and $L\in\R^{N_0}$ be such that $A_0 + \mathfrak{B}_0 K$ and $A_0 - L C_0$ are Hurwitz with eigenvalues that have a real part strictly less than $-\delta<0$. For a given $N \geq N_0 +1$, assume that there exist $P \succ 0$ and $\alpha,\beta,\gamma > 0$ such that 
\begin{equation}\label{eq: thm2 - constraints}
\Theta_1 \preceq 0 ,\quad \Theta_2 \leq 0 ,\quad \Theta_3 \geq 0
\end{equation}
where $\Theta_1$ is defined by (\ref{eq: def theta1 linear}) and 
\begin{align*}
\Theta_2 & = 2\gamma\left\{ -\lambda_{N+1}+q_c + \delta + \frac{1}{\alpha} \right\} + \beta M_\phi \lambda_{N+1}^{3/4} \\
\Theta_3 & = 2\gamma - \frac{\beta M_\phi}{\lambda_{N+1}^{1/4}}
\end{align*}
with $M_\phi = \sum_{n \geq N+1} \frac{\vert \phi_n(0) \vert^2}{\lambda_n^{3/4}} < +\infty$. Then there exists a constant $M > 0$ such that for any initial conditions $z_0 \in H^2(0,1)$ and $\hat{z}_n(0)\in\R$ such that $\cos(\theta_1) z_0(0) - \sin(\theta_1) z_0'(0) = 0$ and $\cos(\theta_2) z_0(1) + \sin(\theta_2) z_0'(1) = K \hat{Z}^{N_0}(0)$, the trajectories of the closed-loop system composed of the plant (\ref{eq: PDE linear}) and the controller (\ref{eq: controller}) satisfy 
$$\Vert z(t,\cdot) \Vert_{L^2}^2 + \sum_{n = 1}^{N} \hat{z}_n(t)^2 \leq M e^{-2\delta t} \left( \Vert z_0 \Vert_{L^2}^2 + \sum_{n = 1}^{N} \hat{z}_n(0)^2 \right)$$ 
for all $t \geq 0$. Moreover, the constraints (\ref{eq: thm2 - constraints}) are always feasible for $N$ selected large enough.
\end{thm}

\begin{pf}
Let the Lyapunov function candidate
\begin{equation}\label{eq: Lyapunov functional - L2 norm}
V(X,w) = X^\top P X + \gamma \sum_{n \geq N+1} \langle w , \phi_n \rangle^2
\end{equation}
for $X \in\R^{2N}$ and $w \in L^2(0,1)$. Proceeding as in the proof of Theorem~\ref{thm1} while replacing the estimate of $\zeta$ by $\zeta^2 \leq M_\phi \sum_{n \geq N+1} \lambda_n^{3/4} w_n^2$, we infer that
$$\dot{V} + 2 \delta V \leq \tilde{X}^\top \Theta_1 \tilde{X} + \sum_{n\geq N+1} \Gamma_n w_n^2$$
where $\tilde{X} = \mathrm{col}\left( X , \zeta \right)$ and $\Gamma_n = 2\gamma\left\{ -\lambda_{n} + q_c + \delta + \frac{1}{\alpha} \right\} + \beta M_\phi \lambda_n^{3/4}$. For $n \geq N+1$ we have $\lambda_n^{3/4} = \lambda_n/\lambda_n^{1/4} \leq \lambda_n/\lambda_{N+1}^{1/4}$ hence 
\begin{align*}
\Gamma_n & \leq - \Theta_3 \lambda_n + 2\gamma \left\{ q_c + \delta + \frac{1}{\alpha} \right\} \\
& \leq - \Theta_3 \lambda_{N+1} + 2\gamma \left\{ q_c + \delta + \frac{1}{\alpha} \right\} = \Theta_2 \leq 0
\end{align*}
where we have used that $\Theta_3 \geq 0$. Combining this result with $\Theta_1 \preceq 0$, we obtain that $\dot{V}+2\delta V \leq 0$ and which implies the claimed stability estimate.

To show that the constraints (\ref{eq: thm2 - constraints}) are always feasible for $N \geq N_0 +1$ selected large enough, we proceed as in the proof of theorem~\ref{thm1} while setting $\alpha = 1$, $\beta = N^{1/8}$, and $\gamma = 1/N^{1/4}$.
\end{pf}

\section{Design in the presence of a sector nonlinearity}\label{sec: Design in the presence of a sector nonlinearity}

Consider now the reaction-diffusion system presenting a sector nonlinearity in the control input, described by
\begin{subequations}\label{eq: PDE nonlinear}
\begin{align}
& z_t(t,x) = (p(x) z_x(t,x))_x - \tilde{q}(x) z(t,x) \\
& \cos(\theta_1) z(t,0) - \sin(\theta_1) z_x(t,0) = 0 \\
& \cos(\theta_2) z(t,1) + \sin(\theta_2) z_x(t,1) = u_\varphi(t) \triangleq \varphi(u(t)) \\
& y (t) = z(t,0) \\
& z(0,x) = z_0(x)
\end{align}
\end{subequations}
for $t > 0$ and $x \in (0,1)$ where $\theta_1 \in (0,\pi/2]$, $\theta_2 \in [0,\pi/2]$, $p \in\mathcal{C}^2([0,1])$ with $p > 0$ and $\tilde{q} \in\mathcal{C}^0([0,1])$. As in the previous section, we consider without loss of generality $q \in\mathcal{C}^0([0,1])$ and $q_c \in\R$ such that (\ref{eq: writting of tilde_q}) holds. The mapping $\varphi : \R \rightarrow \R$ is assumed to be a function of class $\mathcal{C}^1$ for which there exist $k_\varphi > 0$ and $\Delta k_\varphi \in (0,k_\varphi)$ so that
\begin{equation}\label{eq: constraints on varphi}
(k_\varphi - \Delta k_\varphi) \vert x \vert 
\leq \mathrm{sign}(x) \varphi(x) \leq
(k_\varphi + \Delta k_\varphi) \vert x \vert
\end{equation}
for all $x \in\R$. We also assume that $\varphi'$ is locally Lipschitz continuous\footnote{This Lipchitz continuity assumption is only introduced for ensuring the well-posedness of the closed-loop system trajectories but is not explicitly used in the stability analysis.} and $\Vert \varphi' \Vert_{L^\infty} < + \infty$. The objective is to design a finite-dimensional controller and to derive a set of sufficient conditions on the size $\Delta k_\varphi > 0$ of the sector condition (\ref{eq: constraints on varphi}) ensuring the exponential stabilization of (\ref{eq: PDE nonlinear}). 

\begin{rem}\label{rem: nonlinearity}
If one further assume that the mapping $\varphi$ is one to one, then one could merely set $u(t)=\varphi^{-1}(w(t))$ to obtain $u_\varphi(t) = w(t)$. In this particular case, the control design can be easily performed by directly applying any available method for the boundary output feedback stabilization of linear reaction-diffusion equations (such as backstepping, spectral methods, or others). Such an approach fails in general as soon as the mapping $\varphi$ is non-injective, forcing to directly deal with the nonlinear input $u_\varphi(t)=\varphi(u(t))$. In this more general setting, the traditional approach consisting in introducing the time-derivative of the command as an auxiliary input for control design using spectral methods is neither a viable approach. Indeed, defining $v_\varphi = \dot{u}_\varphi = \varphi'(u) \dot{u}$, one could design the auxiliary control input $v_\varphi$. However, the actual control input to be applied remains $u$, which satisfies the dynamics $\varphi'(u(t)) \dot{u}(t) = v_\varphi(t)$ for all $t \geq 0$. Hence, knowing $v_\varphi$, one needs to explicitly compute $u$. This rises two main issues. First, this ODE may not admit solutions defined for all $t \geq 0$. This issue arises, in particular, when $\varphi'$ vanishes at certain isolated points or possibly on certain intervals. Second, provided the well-posedness of this ODE, the stability of this $u$-dynamics is not guaranteed because not assessed when designing $v_\varphi$.
\end{rem}

\subsection{Spectral reduction}
Introducing the change of variable
\begin{equation}\label{eq: change of variable - nonlinear}
w(t,x) = z(t,x) - \frac{x^2}{\cos(\theta_2) + 2 \sin(\theta_2)} u_\varphi(t) 
\end{equation}
we infer that
\begin{subequations}\label{eq: PDE nonlinear - homogeneous}
\begin{align}
& v_\varphi(t) = \dot{u}_\varphi(t) = \varphi'(u(t)) \dot{u}(t) \\
& w_t(t,x) = (p(x) w_x(t,x))_x - \tilde{q}(x) w(t,x) \\
& \phantom{w_t(t,x) =} \;  + a(x) u_\varphi(t) + b(x) v_\varphi(t) \\ 
& \cos(\theta_1) w(t,0) - \sin(\theta_1) w_x(t,0) = 0 \\
& \cos(\theta_2) w(t,1) + \sin(\theta_2) w_x(t,1) = 0 \\
& y(t) = w(t,0) \\
& w(0,x) = w_0(x)
\end{align}
\end{subequations}
where $a,b$ are defined as in the previous section and $w_0(x) = z_0(x) - \frac{x^2}{\cos(\theta_2) + 2 \sin(\theta_2)} u_\varphi(0)$. With the coefficients of projection defined in the previous section, we have
\begin{equation}\label{eq: link z_n and w_n - nonlinear}
w_n(t) = z_n(t) + b_n u_\varphi(t), \quad n \geq 1 .
\end{equation}
The projection of (\ref{eq: PDE nonlinear}) into $(\phi_n)_{n \geq 1}$ gives
\begin{equation}\label{eq: dynamics z_n - nonlinear}
\dot{z}_n(t) = (-\lambda_n + q_c) z_n(t) + \beta_n u_\varphi(t)
\end{equation}
with $\beta_n$ defined by (\ref{eq: def beta_n}) while the projection of (\ref{eq: PDE nonlinear - homogeneous}) reads
\begin{subequations}\label{eq: dynamics w_n - nonlinear}
\begin{align}
\dot{u}_\varphi & = v_\varphi = \varphi'(u) \dot{u} \\
\dot{w}_n(t) & = (-\lambda_n + q_c) w_n(t) + a_n u_\varphi(t) + b_n v_\varphi(t) , \; n \geq 1 \\
y(t) & = \sum_{n \geq 1} w_n(t) \phi_n(0) \label{eq: dynamics w_n - y - nonlinear}
\end{align}
\end{subequations}

\begin{rem}
Representation (\ref{eq: PDE nonlinear - homogeneous}) cannot be used for control design with $v_\varphi$ selected as an auxiliary input signal. This is because, as discussed in Remark~\ref{rem: nonlinearity}, $v_\varphi = \varphi'(u)\dot{u}$ where $u$ remains the actual to-be-implemented input of the plant (\ref{eq: PDE nonlinear}). Hence the approach proposed in~\cite{lhachemi2020finite} is inapplicable in the presence of the input nonlinearity $\varphi$. We solve this problem by adopting the approach reported in the previous section, namely by performing the control design on (\ref{eq: PDE nonlinear}) while carrying out the Lyapunov-based stability analysis using (\ref{eq: PDE nonlinear - homogeneous}).
\end{rem}

\subsection{Control strategy}
Let $\delta > 0$ and $N_0 \geq 1$ be such that $-\lambda_n + q_c < - \delta < 0$ for all $n \geq N_0 + 1$. We consider the following observer-based control strategy: 
\begin{subequations}\label{eq: controller - nonlinear}
\begin{align}
\hat{w}_n(t) & = \hat{z}_n(t) + b_n u_\varphi (t) \label{eq: controller 1 - nonlinear} \\
\dot{\hat{z}}_n(t) & = (-\lambda_n+q_c) \hat{z}_n(t) + \beta_n u_\varphi(t) \label{eq: controller 2 - nonlinear} \\
& \phantom{=}\; - l_n \left\{ \sum_{k = 1}^N \hat{w}_k(t) \phi_k(0) - y(t) \right\}  ,\; 1 \leq n \leq N_0 \nonumber \\
\dot{\hat{z}}_n(t) & = (-\lambda_n+q_c) \hat{z}_n(t) + \beta_n u_\varphi(t) ,\; N_0+1 \leq n \leq N \label{eq: controller 3 - nonlinear} \\
u(t) & = \sum_{k = 1}^{N_0} k_k \hat{z}_k(t) \label{eq: controller 4 - nonlinear}
\end{align}
\end{subequations}
Here $l_n \in\R$ and $k_k \in\R$ are the observer and feedback gains, respectively. 

\begin{rem}\label{rem: well-posedness}
We denote $\hat{z}(t) \in\R^N$ the state of the observer. Under the above mentioned assumption for the sector nonlinearity $\varphi$, the well-posedness of the closed-loop system composed on (\ref{eq: PDE nonlinear - homogeneous}) and (\ref{eq: controller - nonlinear}) in terms of classical solutions for initial conditions $w_0 \in D(\mathcal{A})$ and $\hat{z}(0) \in\R^N$, namely $(w,\hat{z}) \in \mathcal{C}^0([0,T);L^2(0,1) \times \R^N) \cap \mathcal{C}^1((0,T);L^2(0,1) \times \R^N)$, defined on a maximal interval of existence $[0,T)$ with either $T > 0$ or $T = + \infty$, is a direct consequence of~\cite[Thm.~6.3.1]{pazy2012semigroups}. Moreover, $w(t,\cdot) \in D(\mathcal{A})$ for all $t > 0$ and, from the proof of~\cite[Thm.~6.3.1]{pazy2012semigroups}, we have $\mathcal{A} w \in \mathcal{C}^0((0,T);L^2(0,1))$ and $\mathcal{A}^{1/2}w \in \mathcal{C}^0([0,T);L^2(0,1))$. Finally, from a similar argument that the one stated as a preliminary remark of the proof of~\cite[Thm.~6.3.3]{pazy2012semigroups}, if $T < +\infty$ then $\Vert \mathcal{A}^{1/2} w(t,\cdot) \Vert_{L^2}^2 + \Vert \hat{z}(t) \Vert^2 = \sum_{n \geq 1} \lambda_n w_n(t)^2 + \Vert \hat{z}(t) \Vert^2$ is unbounded on $[0,T)$.
\end{rem}

\subsection{Model for stability analysis}

We define the mapping $\psi : \R\rightarrow\R$ by
\begin{equation}\label{eq: def psi}
\psi(x) = \varphi(x) - k_\varphi x .
\end{equation}
Adopting the definitions and the approach of Subsection~\ref{subsec: Model for stability analysis}, we infer from (\ref{eq: controller - nonlinear}) that 
\begin{subequations}\label{eq: truncated model nonlinear - 4 ODEs}
\begin{align}
u & = K \hat{Z}^{N_0} \label{eq: truncated model nonlinear - 4 ODEs - 0} \\
\dot{\hat{Z}}^{N_0} & = (A_0 + k_\varphi \mathfrak{B}_0 K) \hat{Z}^{N_0} + LC_0 E^{N_0} \label{eq: truncated model nonlinear - 4 ODEs - 1} \\
& \phantom{=}\; + L\tilde{C}_1 \tilde{E}^{N-N_0} + L \zeta + \mathfrak{B}_0 \psi(K\hat{Z}^{N_0}) \nonumber \\
\dot{E}^{N_0} & = ( A_0 - L C_0 ) E^{N_0} - L \tilde{C}_1 \tilde{E}^{N-N_0} - L \zeta \\
\dot{\tilde{Z}}^{N-N_0} & = A_1 \tilde{Z}^{N-N_0} + k_\varphi \tilde{\mathfrak{B}}_1 K \hat{Z}^{N_0} + \tilde{\mathfrak{B}}_1 \psi(K\hat{Z}^{N_0}) \\
\dot{\tilde{E}}^{N-N_0} & = A_1 \tilde{E}^{N-N_0}
\end{align}
\end{subequations}
Hence, with $X$ defined by (\ref{eq: truncated model - def X}), $u$ can still be expressed by $u = \tilde{K} X$ and we have that 
\begin{equation}\label{eq: truncated model - nonlinear}
\dot{X} = F X + \mathcal{L} \zeta + \mathcal{L}_\psi \psi(K\hat{Z}^{N_0}) 
\end{equation}
where
\begin{equation*}
F =
\begin{bmatrix}
A_0 + k_\varphi \mathfrak{B}_0 K & LC_0 & 0 & L\tilde{C}_1 \\
0 & A_0 - L C_0 & 0 & - L\tilde{C}_1 \\
k_\varphi\tilde{\mathfrak{B}}_1 K & 0 & A_1 & 0 \\
0 & 0 & 0 & A_1
\end{bmatrix} ,
\end{equation*}
$\mathcal{L} = \mathrm{col}\left( L , -L , 0 , 0 \right)$, and $\mathcal{L}_\psi = \mathrm{col}\left( \mathfrak{B}_0 , 0 , \tilde{\mathfrak{B}}_1 , 0 \right)$. With $\tilde{X} = \mathrm{col}\left( X , \zeta , \psi(K\hat{Z}^{N_0}) \right)$ 
and based on (\ref{eq: truncated model nonlinear - 4 ODEs - 0}-\ref{eq: truncated model nonlinear - 4 ODEs - 1}), we infer that 
\begin{equation}\label{eq: derivative v of command input u - nonlinear}
v_\varphi = \dot{u}_\varphi = \varphi'(K \hat{Z}^{N_0}) K \dot{\hat{Z}}^{N_0} = \varphi'(K \hat{Z}^{N_0}) E \tilde{X}
\end{equation}
where $E = K \begin{bmatrix} A_0 + k_\varphi \mathfrak{B}_0 K & LC_0 & 0 & L\tilde{C}_1 & L & \mathfrak{B}_0 \end{bmatrix}$.

\subsection{Main stability results}

\begin{thm}\label{thm3}
Let $\theta_1 \in (0,\pi/2]$, $\theta_2 \in [0,\pi/2]$, $p \in\mathcal{C}^2([0,1])$ with $p > 0$, and $\tilde{q} \in\mathcal{C}^0([0,1])$. Let $k_\varphi > 0$, $\Delta k_\varphi \in (0,k_\varphi)$, and $M_\varphi > 0$. Let $q \in\mathcal{C}^0([0,1])$ and $q_c \in\R$ be such that (\ref{eq: writting of tilde_q}) holds. Let $\delta > 0$ and $N_0 \geq 1$ be such that $-\lambda_n + q_c < - \delta$ for all $n \geq N_0 + 1$. Let $K\in\R^{1 \times N_0}$ and $L\in\R^{N_0}$ be such that $A_0 + k_\varphi \mathfrak{B}_0 K$ and $A_0 - L C_0$ are Hurwitz with eigenvalues that have a real part strictly less than $-\delta<0$. For a given $N \geq N_0 +1$, assume that there exist $P \succ 0$, $\alpha>3/2$, and $\beta,\gamma,\tau > 0$ such that 
\begin{equation}\label{eq: thm3 - constraints}
\Theta_1 \preceq 0 ,\quad \Theta_2 \leq 0
\end{equation}
where
\begin{align*}
\Theta_1 & = \begin{bmatrix} \Theta_{1,1} & P\mathcal{L} & P\mathcal{L}_\psi \\ \mathcal{L}^\top P & - \beta & 0 \\ \mathcal{L}_\psi^\top P & 0 & \alpha\gamma \Vert \mathcal{R}_N a \Vert_{L^2}^2 - \tau \end{bmatrix} \\
& \phantom{=}\; + \alpha\gamma \Vert \mathcal{R}_N b \Vert_{L^2}^2 M_\varphi^2 E^\top E \\
\Theta_{1,1} & = F^\top P + P F + 2 \delta P \\
& \phantom{=}\; + \left\{ \alpha\gamma k_\varphi^2 \Vert \mathcal{R}_N a \Vert_{L^2}^2 + \tau \Delta k_\varphi^2 \right\} \tilde{K}^\top \tilde{K} \\
\Theta_2 & = 2\gamma\left\{ - \left( 1 - \frac{3}{2\alpha} \right) \lambda_{N+1}+q_c + \delta \right\} + \beta M_\phi
\end{align*}
and with $M_\phi = \sum_{n \geq N+1} \frac{\vert \phi_n(0) \vert^2}{\lambda_n} < +\infty$. Then there exists a constant $M > 0$ such that for any $\varphi \in\mathcal{C}^1(\R)$ such that (\ref{eq: constraints on varphi}) holds with $\varphi'$ locally Lipschitz continuous and $\Vert \varphi' \Vert_{L^\infty} \leq M_\varphi$, and for any initial conditions $z_0 \in H^2(0,1)$ and $\hat{z}_n(0)\in\R$ such that $\cos(\theta_1) z_0(0) - \sin(\theta_1) z_0'(0) = 0$ and $\cos(\theta_2) z_0(1) + \sin(\theta_2) z_0'(1) = \varphi(K \hat{Z}^{N_0}(0))$, the trajectories of the closed-loop system composed of the plant (\ref{eq: PDE nonlinear}) and the controller (\ref{eq: controller - nonlinear}) satisfy 
$$\Vert z(t,\cdot) \Vert_{H^1}^2 + \sum_{n = 1}^{N} \hat{z}_n(t)^2 \leq M e^{-2\delta t} \left( \Vert z_0 \Vert_{H^1}^2 + \sum_{n = 1}^{N} \hat{z}_n(0)^2 \right)$$ 
for all $t \geq 0$. Moreover, for any given $k_\varphi,M_\varphi > 0$, there exists $\Delta k_\varphi \in (0,k_\varphi)$ such that the constraints (\ref{eq: thm3 - constraints}) are always feasible when selecting $N$ to be large enough.
\end{thm}

\begin{pf}
Considering the Lyapunov function candidate defined by (\ref{eq: Lyapunov functional - H1 norm}), the computation of its time derivative along the system trajectories (\ref{eq: dynamics w_n - nonlinear}) and (\ref{eq: truncated model - nonlinear}) gives
\begin{align*}
\dot{V} & =
\tilde{X}^\top
\begin{bmatrix} F^\top P + P F & P \mathcal{L} & P \mathcal{L}_\psi \\ \mathcal{L}^\top P & 0 & 0 \\ \mathcal{L}_\psi^\top P & 0 & 0 \end{bmatrix}
\tilde{X} \\
& \phantom{=}\; + 2 \gamma \sum_{n \geq N+1} \lambda_n \left\{ (-\lambda_n + q_c) w_n + a_n u_\varphi + b_n v_\varphi \right\} w_n .
\end{align*}
where $\tilde{X} = \mathrm{col}\left( X , \zeta , \psi(K\hat{Z}^{N_0}) \right)$. Since $u_\varphi = \varphi(K\hat{Z}^{N_0}) = k_\varphi \tilde{K}X + \psi(K\hat{Z}^{N_0})$, using Young inequality, we infer for any $\alpha > 0$ that
\begin{align*}
& 2\sum_{n \geq N+1} \lambda_n a_n u_\varphi w_n \leq \frac{2}{\alpha} \sum_{n \geq N+1} \lambda_n^2 w_n^2 \\
& \quad\qquad + \alpha \Vert \mathcal{R}_N a \Vert_{L^2}^2 \left\{ k_\varphi^2 (\tilde{K}X)^2 + \psi(K\hat{Z}^{N_0})^2 \right\} , \\
& 2\sum_{n \geq N+1} \lambda_n b_n v_\varphi w_n \leq \frac{1}{\alpha} \sum_{n \geq N+1} \lambda_n^2 w_n^2 + \alpha \Vert \mathcal{R}_N b\Vert_{L^2}^2 v_\varphi^2 
\end{align*}
where, using (\ref{eq: derivative v of command input u - nonlinear}), $v_\varphi^2 \leq M_\varphi^2 \tilde{X}^\top E^\top E \tilde{X}$. Moreover, since $\zeta = \sum_{n \geq N+1} w_n \phi_n(0)$, Cauchy-Schwarz inequality gives $\zeta^2 \leq M_\phi \sum_{n \geq N+1} \lambda_n w_n^2$. Combining the above estimates, we have
$\dot{V} + 2 \delta V \leq \tilde{X}^\top \Theta_{1,\tau=0} \tilde{X} + \sum_{n\geq N+1}\lambda_n \Gamma_n w_n^2$
where $\Theta_{1,\tau=0}$ is obtained from $\Theta_1$ by setting $\tau = 0$ and $\Gamma_n = 2\gamma\left\{ - \left( 1 - \frac{3}{2\alpha} \right) \lambda_{n}+q_c + \delta \right\} + \beta M_\phi$. We now need to take advantage of the sector condition (\ref{eq: constraints on varphi}) satisfied by $\varphi$. More precisely, (\ref{eq: constraints on varphi}) implies that $( \varphi(x) - (k_\varphi + \Delta k_\varphi) x ) ( \varphi(x) - (k_\varphi - \Delta k_\varphi) x ) \leq 0$ for all $x \in\R$. Using (\ref{eq: def psi}), this is equivalent to $\psi(x)^2 - \Delta k_\varphi^2 x^2 \leq 0$ for all $x \in\R$. The use of this sector condition applied at $x = K \hat{Z}^{N_0} = \tilde{K} X$ into the above estimate of $\dot{V} + 2 \delta V$ implies that 
$$\dot{V} + 2 \delta V \leq \tilde{X}^\top \Theta_{1} \tilde{X} + \sum_{n\geq N+1}\lambda_n \Gamma_n w_n^2 .$$
Using $\alpha > 3/2$, we have $\Gamma_n \leq \Theta_2 \leq 0$ for all $n \geq N+1$. Since $\Theta_1 \preceq 0$, we infer that $\dot{V} + 2 \delta V \leq 0$. From (\ref{eq: Lyapunov functional - H1 norm}) we infer that $\Vert \mathcal{A}^{1/2} w(t,\cdot) \Vert_{L^2}^2 = \sum_{n \geq 1} \lambda_n w_n(t)^2$ and $\Vert \hat{z}(t) \Vert$ are bounded on the maximal interval of existence of the system trajectories. Hence the trajectories are well-defined for all $t \geq 0$ (see end of Remark~\ref{rem: well-posedness}) and we obtain the claimed stability estimate.

It remains to show for any given $k_\varphi,M_\varphi > 0$ the existence of some $\Delta k_\varphi \in (0,k_\varphi)$ so that the constraints $\Theta_1 \preceq 0$ and $\Theta_2 \leq 0$ are always feasible when taking $N \geq N_0 + 1$ large enough. Proceeding as in the proof of Theorem~\ref{thm1}, the application of Lemma~\ref{lem: useful lemma} reported in appendix to the matrix $F + \delta I$ ensures that the solution $P \succ 0$ to the Lyapunov equation $F^\top P + P F + 2 \delta P = -I$ is such that $\Vert P \Vert = O(1)$ as $N \rightarrow + \infty$. This ensures the existence of a constant $M_P > 0$ such that $\Vert P \Vert \leq M_P$ for all $N \geq N_0 + 1$. There also exists a constant $M_\psi > 0$ such that $\Vert \mathcal{L}_\psi \Vert \leq M_\psi$ for all $N \geq N_0 + 1$. This allows us to define $\tau = 1 + 4 M_P^2 M_\psi^2 + \Vert a \Vert_{L^2}^2$ and $\Delta k_\varphi = \min\left( \frac{1}{1+2 \Vert K \Vert \sqrt{\tau}},\frac{k_\varphi}{2}\right)$ which are constants independent of $N$ and $\varphi$. We also set $\alpha = 2$, $\beta = \sqrt{N}$, and $\gamma = 1/N$. Since $\Vert \tilde{K} \Vert = \Vert K \Vert$ and $\Vert \mathcal{L} \Vert = \sqrt{2} \Vert L \Vert$ are constants independent of $N$, the use of Schur complement implies that
\begin{equation*}
\Xi_1 = \begin{bmatrix}
-I + \alpha\gamma k_\varphi^2 \Vert \mathcal{R}_N a \Vert_{L^2}^2 \tilde{K}^\top \tilde{K} & P \mathcal{L} \\
\mathcal{L}^\top P & - \beta
\end{bmatrix}
\preceq - \frac{3}{4} I .
\end{equation*}
for $N$ large enough. Since $\tau - \alpha\gamma \Vert \mathcal{R}_N a \Vert_{L^2}^2 - 1/2 > 0$ for $N \geq 2$,  we infer from the Schur complement that
\begin{equation*}
\Xi_2 = \begin{bmatrix} \Xi_1 & \Psi \\ \Psi^\top & \alpha\gamma \Vert \mathcal{R}_N a \Vert_{L^2}^2 - \tau \end{bmatrix} \preceq -\frac{1}{2} I , 
\end{equation*}
where $\Psi = \mathrm{col}(P\mathcal{L}_\psi,0)$,
if and only if 
$$\Xi_1 + \frac{1}{2} I + \frac{1}{\tau - \alpha\gamma \Vert \mathcal{R}_N a \Vert_{L^2}^2 - 1/2} \Psi \Psi^\top  \preceq 0 .$$ 
A sufficient condition ensuring that this latter inequality is satisfied is provided by $4 \Vert P \Vert^2 \Vert \mathcal{L}_\psi \Vert^2 \leq 4 M_P^2 M_\psi^2 \leq \tau - \alpha\gamma \Vert \mathcal{R}_N a \Vert_{L^2}^2 - \frac{1}{2}$, which is true for $N \geq 2$ based on the definitions of $\tau,\alpha,\gamma$. Noting now from (\ref{eq: derivative v of command input u - nonlinear}) that $\Vert E \Vert = O(1)$ as $N \rightarrow + \infty$, we obtain that $\alpha\gamma \Vert \mathcal{R}_N b \Vert_{L^2}^2 M_\varphi^2 E^\top E \preceq \frac{1}{4}I$ for $N$ large enough. Putting together the above estimates, we infer that $\Theta_1 \preceq -\frac{1}{4} I + \tau \Delta k_\varphi^2 \Vert K \Vert^2 I$ for $N$ large enough. Since $\Delta k_\varphi \leq 1 / (1 + 2 \Vert K \Vert \sqrt{\tau})$, we have $\Theta_1 \preceq 0$ for $N$ large enough. Finally, recalling that $\alpha = 2$, $\beta = \sqrt{N}$, and $\gamma = 1/N$, we also note that $\Theta_2 \leq 0$ for $N$ large enough. In conclusion, we have found a $\Delta k_\varphi \in (0,k_\varphi)$ such that the constraints (\ref{eq: thm3 - constraints}) are feasible when selecting $N \geq N_0 + 1$ to be large enough. This completes the proof
\end{pf}

\begin{rem}
A similar approach that the one reported in Remark~\ref{rem: obtention of LMIs} can be employed to recast the constraints (\ref{eq: thm3 - constraints}) into a LMI formulation for any given $N \geq N_0 + 1$.
\end{rem}

\begin{rem}
For any arbitrary choice of the feedback and observer gains $K \in\R^{1 \times N_0}$ and $L \in\R^{N_0}$ such that the matrices $A_0 + k_\varphi \mathfrak{B}_0 K$ and $A_0 - L C_0$ are Hurwitz, Theorem~\ref{thm3} ensures the existence and allows the computation of a $\Delta k_\varphi> 0$. Hence, the maximum value obtained for $\Delta k_\varphi> 0$ depends on the specific realizations of the gains $K$ and $L$. Note however that the problem of maximizing the value of $\Delta k_\varphi> 0$ by tuning the gains $K$ and $L$ is difficult because highly nonlinear. For example, $\Theta_1$ presents through the term $E^\top E$ a biquadratic dependency on the feedback gain $K$.
\end{rem}

As a corollary of Theorem~\ref{thm3}, we have the following $L^2$ version for the stability of the closed-loop system.

\begin{cor}\label{thm4}
In addition of all the assumptions of Theorem~\ref{thm3}, assume further that $N \geq N_0 +1$ is selected such that there exist $P' \succ 0$, and $\alpha',\beta',\gamma',\tau' > 0$ so that 
\begin{equation}\label{eq: thm4 - constraints}
\Theta_1' \preceq 0 ,\quad \Theta_2' \leq 0 ,\quad \Theta_3' \geq 0
\end{equation}
where
\begin{align*}
\Theta_1' & = \begin{bmatrix} \Theta_{1,1}' & P'\mathcal{L} & P'\mathcal{L}_\psi \\ \mathcal{L}^\top P' & - \beta' & 0 \\ \mathcal{L}_\psi^\top P' & 0 & \alpha'\gamma' \Vert \mathcal{R}_N a \Vert_{L^2}^2 - \tau' \end{bmatrix} \\
& \phantom{=}\; + \alpha'\gamma' \Vert \mathcal{R}_N b \Vert_{L^2}^2 M_\varphi^2 E^\top E \\
\Theta_{1,1}' & = F^\top P' + P' F + 2 \delta P' \\
& \phantom{=}\; + \left\{ \alpha'\gamma' k_\varphi^2 \Vert \mathcal{R}_N a \Vert_{L^2}^2 + \tau' \Delta k_\varphi^2 \right\} \tilde{K}^\top \tilde{K} \\
\Theta_2' & = 2\gamma'\left\{ -\lambda_{N+1}+q_c + \delta + \frac{3}{2\alpha'} \right\} + \beta' M_\phi' \lambda_{N+1}^{3/4} \\
\Theta_3' & = 2\gamma' - \frac{\beta' M_\phi'}{\lambda_{N+1}^{1/4}}
\end{align*}
and with $M_\phi' = \sum_{n \geq N+1} \frac{\vert \phi_n(0) \vert^2}{\lambda_n^{3/4}} < +\infty$. Then there exists a constant $M' > 0$ such that, under the same assumptions for $\varphi$ and the initial conditions that the ones of Theorem~\ref{thm3}, the trajectories of the closed-loop system composed of the plant (\ref{eq: PDE nonlinear}) and the controller (\ref{eq: controller - nonlinear}) satisfy 
$$\Vert z(t,\cdot) \Vert_{L^2}^2 + \sum_{n = 1}^{N} \hat{z}_n(t)^2 \leq M' e^{-2\delta t} \left( \Vert z_0 \Vert_{L^2}^2 + \sum_{n = 1}^{N} \hat{z}_n(0)^2 \right)$$ 
for all $t \geq 0$. Moreover, for any given $k_\varphi,M_\varphi > 0$, there exists $\Delta k_\varphi \in (0,k_\varphi)$ such that the both constraints (\ref{eq: thm3 - constraints}) and (\ref{eq: thm4 - constraints}) are always feasible when selecting $N$ to be large enough.
\end{cor}

\begin{pf}
From Theorem~\ref{thm3}, we infer the existence of solutions defined for all $t \geq 0$. Now the proof is based on the Lyapunov functional (\ref{eq: Lyapunov functional - L2 norm}) and relies on similar arguments that the ones reported for Theorems~\ref{thm2} and~\ref{thm3}.
\end{pf}

\begin{rem}
The statement of Corollary~\ref{thm4} requires the constraints (\ref{eq: thm4 - constraints}) of Theorem~\ref{thm3}. This is because, even if (\ref{eq: thm4 - constraints}) is sufficient by itself to ensure the exponential decay in $L^2$ norm, it fails to ensure that $\Vert \mathcal{A}^{1/2} w(t,\cdot) \Vert_{L^2}^2 = \sum_{n \geq 1} \lambda_n w_n(t)^2$ remains bounded on any time intervals of finite length. This latter point is required to apply the argument at the end of Remark~\ref{rem: well-posedness} to ensure that the system trajectories are well defined for all $t \geq 0$.
\end{rem}

\section{Numerical illustration}\label{sec: numerical illustration}

\begin{figure}
\center
\includegraphics[width=3.5in]{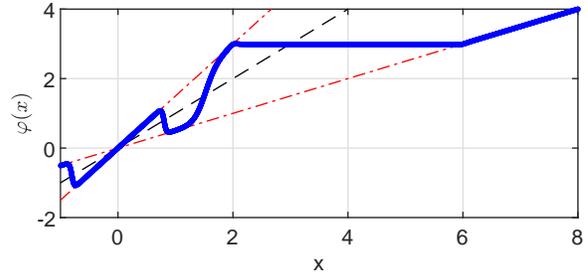}
\caption{Input sector nonlinear $\varphi$ with $\varphi(x) \sim 0.5 x$ for $x < -6$ and $x > 6$}
\label{fig: sector nonlinearity}
\end{figure}

Let the unstable reaction-diffusion equation described by (\ref{eq: PDE nonlinear}) with $\theta_1 = \pi/2$, $\theta_2 = 0$, $p = 1$ and $\tilde{q} = - 3$. We consider the sector nonlinearity $\varphi$ that takes the form depicted in Fig~\ref{fig: sector nonlinearity} and which satisfies (\ref{eq: constraints on varphi}) for $k_\varphi = 1$ and $\Delta k_\varphi = 0.5$ while $\Vert \varphi' \Vert_{L^\infty} \leq 9.02$. This function is not injective and has a derivative that vanishes both at certain isolated points and on an interval. Hence, this nonlinear mapping checks all the pathological behaviors discussed in Remark~\ref{rem: nonlinearity} that prevent the design of the control strategy using the auxiliary control input $v_\varphi = \dot{u}_\varphi$.

\begin{figure}
     \centering
     	\subfigure[State of the reaction-diffusion system $z(t,x)$]{
		\includegraphics[width=3.5in]{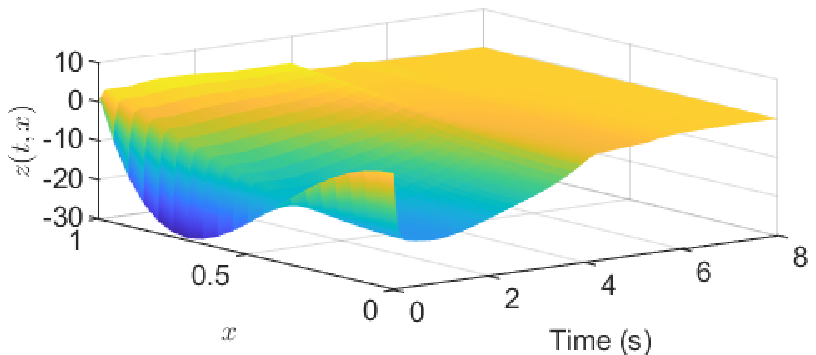}
		}
     	\subfigure[Error of observation $e(t,x)$]{
		\includegraphics[width=3.5in]{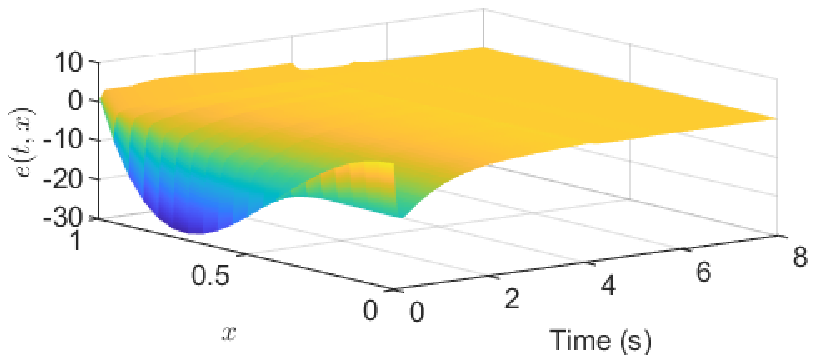}
		}
     	\subfigure[Control input]{
		\includegraphics[width=3.5in]{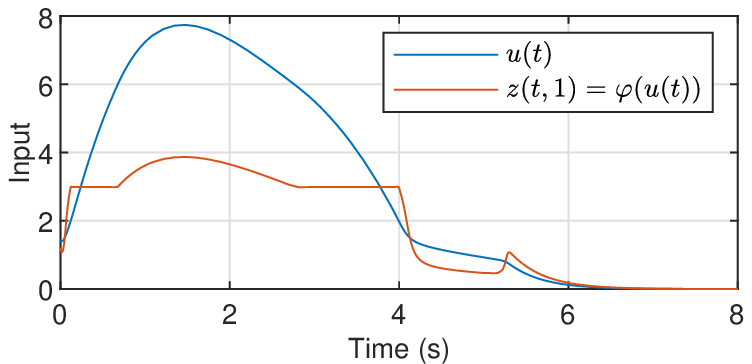}
		\label{fig: sim CL - input}
		}		
     \caption{Time evolution of the closed-loop system}
     \label{fig: sim CL}
\end{figure}

With $N_0 = 1$, the feedback gain $K = -0.8250$, and the observer gain $L = 1.2958$, the sufficient conditions (\ref{eq: thm3 - constraints}) of Theorem~\ref{thm3} are found feasible for $\delta = 0.3$ when using an observer of dimension $N = 3$, ensuring the exponential stability of the closed-loop system in $H^1$-norm. Simulation results are presented in Fig~\ref{fig: sim CL}. It can be observed that the both state of the PDE and observation error converge exponentially to zero despite the strong impact of the sector nonlinearity applying on the control input as shown in Fig.~\ref{fig: sim CL - input}. This is compliant with the theoretical predictions of Theorem~\ref{thm3}. For the same system, the conditions of Corollary~\ref{thm4} that ensure the stability of the system in $L^2$-norm appear to be more stringent from a numerical perspective. For instance, the sufficient conditions (\ref{eq: thm3 - constraints}) and (\ref{eq: thm4 - constraints}) of Corollary~\ref{thm4} are found feasible for $\delta = 0.3$ when using an observer of dimension $N = 16$.

It can be observed in simulation that the reported control strategy cannot achieve the stabilization of the plant for arbitrarily large values of the size $\Delta k _\varphi$ of the sector non linearity (\ref{eq: constraints on varphi}). Indeed, in the setting of the previous paragraph with a nonlinearity similar to the one of Fig~\ref{fig: sector nonlinearity} but rescaled with a size of the sector condition increased from $\Delta k _\varphi = 0.5$ to $\Delta k _\varphi \approx 0.72$, we obtain in simulation divergent closed-loop system trajectories. Moreover, one can expect that the value of $\Delta k_\varphi$ decreases with the level of instability of the open-loop plant (characterized by its growth bound), i.e, when decreasing the value of $\tilde{q}$. This can be observed numerically by considering for example the constraints of Theorem~3 with $k_\varphi = 1$ and $\delta = 0.3$ while placing the poles of both $A_0 + k_\varphi \mathfrak{B}_0 K$ and $A_0 - L C_0$ at $-1.3$. In this case we obtain for a dimension $N = 15$ of the observer: $\Delta k_\varphi = 0.54$ for $\tilde{q}=-3$, $\Delta k_\varphi = 0.24$ for $\tilde{q}=-5$, $\Delta k_\varphi = 0.12$ for $\tilde{q}=-7$, and $\Delta k_\varphi = 0.03$ for $\tilde{q}=-9$.

\section{Conclusion}\label{sec: conclusion}

This paper has investigated the finite-dimensional observer-based stabilization of a reaction-diffusion equation in the presence of a sector nonlinearity in the control input. It is worth noting that even if the method has been presented in the case of a Robin boundary input with parameter $\theta_1 \in(0,\pi/2]$ and $\theta_2 \in[0,\pi/2]$, the approach readily extends to the case $\theta_1 \in(0,\pi)$ and $\theta_2 \in[0,\pi)$ provided $q$ in (\ref{eq: writting of tilde_q}) is selected sufficiently large positive so that (\ref{eq: inner product Af and f}) still holds and by replacing the change of variable (\ref{eq: change of variable}) by $w(t,x) = z(t,x) - \frac{x^\alpha}{\cos(\theta_2) + \alpha \sin(\theta_2)} u(t)$ for any fixed $\alpha > 1$ so that $\cos(\theta_2) + \alpha \sin(\theta_2) \neq 0$. Future research directions may be concern with extensions to other types of nonlinearities.


\bibliographystyle{plain}        
\bibliography{autosam}           




\appendix
\section{Useful Lemma}\label{annex: proof lemma well-posedness}

The following Lemma is borrowed from~\cite{lhachemi2020finite} and generalizes a result presented in~\cite{katz2020constructive}.

\begin{lem}\label{lem: useful lemma}
Let $n,m,N \geq 1$, $M_{11} \in \R^{n \times n}$ and $M_{22} \in \R^{m \times m}$ Hurwitz, $M_{12} \in \R^{n \times m}$, $M_{14}^N \in\R^{n \times N}$, $M_{24}^N \in\R^{m \times N}$, $M_{31}^N \in\R^{N \times n}$, $M_{33}^N,M_{44}^N \in \R^{N \times N}$, and
\begin{equation*}
F^N = \begin{bmatrix}
M_{11} & M_{12} & 0 & M_{14}^N \\
0 & M_{22} & 0 & M_{24}^N \\
M_{31}^N & 0 & M_{33}^N & 0 \\
0 & 0 & 0 & M_{44}^N
\end{bmatrix} .
\end{equation*}
We assume that there exist constants $C_0 , \kappa_0 > 0$ such that $\Vert e^{M_{33}^N t} \Vert \leq C_0 e^{-\kappa_0 t}$ and $\Vert e^{M_{44}^N t} \Vert \leq C_0 e^{-\kappa_0 t}$ for all $t \geq 0$ and all $N \geq 1$. Moreover, we assume that there exists a constant $C_1 > 0$ such that $\Vert M_{14}^N \Vert \leq C_1$, $\Vert M_{24}^N \Vert \leq C_1$, and $\Vert M_{31}^N \Vert \leq C_1$ for all $N \geq 1$. Then there exists a constant $C_2 > 0$ such that, for any $N \geq 1$, there exists a symmetric matrix $P^N \in\R^{n+m+2N}$ with $P^N \succ 0$ such that $(F^N)^\top P^N + P^N F^N = - I$ and $\Vert P^N \Vert \leq C_2$.
\end{lem}

\end{document}